\documentclass[11pt,letter,subeqn]{article}
\usepackage{amsmath, amssymb,amsthm,cite}
\usepackage{graphics,color}
\usepackage{graphicx}
\usepackage{tabulary}
\usepackage{float}
\usepackage{caption, subcaption}
\usepackage{subfiles}
\usepackage{multirow}
\usepackage{chngcntr}
\usepackage{cleveref}
\counterwithin{figure}{section}
\counterwithin{table}{section}

\title{Approximate Symmetries and Approximate Solutions of Some Perturbed ODE Models}

\author{ 
Mahmood R. Tarayrah\vspace{0.5cm}\\
\emph{Department of Mathematics and Statistics, University of Saskatchewan, Saskatoon, Canada}\vspace{0.2cm}\\
Electronic mail: mrt566@usask.ca
}

\def\beq{\begin{equation}}
\def\eeq{\end{equation}}
\def\barr{\begin{array}{ll}}
\def\earr{\end{array}}

\def\sech{\mathop{\hbox{\rm sech}}}

\newtheorem{theorem}{Theorem}[section]

{\theoremstyle{definition}

}

\newcounter{tabnum}

\begin{document}


\maketitle \numberwithin{equation}{section}
\maketitle \numberwithin{remark}{section}
\numberwithin{lemma}{section}
\numberwithin{proposition}{section}

\begin{abstract}

We find Baikov-Gazizov-Ibragimov approximate point symmetries of the second-order Boussinesq ODE, and we find the higher-order approximate symmetries corresponding to the unstable point symmetries (the point symmetries that disappear fron the classification of the BGI approximate point symmetries) of the unperturbed equation. Approximate local symmetries are used to construct a general approximate solution of the Boussinesq ODE. We use approximate integrating factors to find a general approximate solution of the Benjamin-Bona-Mahony ODE reduction.

\end{abstract}

{\bf Keywords:}~ Lie groups, Local symmetries, Approximate symmetries; Ordinary differential equations, Approximate solutions.

\section{Introduction}
A symmetry of a system of differential equations is a transformation that maps solutions of the system to other solutions. Symmetry ideas have been developed over the last century, relating to symmetry reduction and solution of differential equations, integrating factors, conserved quantities and local conservation laws, and more (see, e.g., Refs.~\cite{olver2000applications, bluman2010applications} and references therein).

Perturbed differential equations are differential equations involving a small parameter. This small perturbation disturbs the symmetry group properties of the unperturbed equations. A class of Lie symmetries which are useful in studying the symmetry properties of the perturbed differential equations and/or provide new symmetries for these equations are the approximate Lie symmetries. Several approximate Lie symmetry methods have been developed. The approximate symmetry group method was introduced by Baikov, Gazizov and Ibragimov \cite{baikov1989, baikov1991, baikov1993}, we shall call their method the \textrm{BGI} method. They introduced a one-parameter approximate transformation group where the approximate symmetry generator is expanded in a perturbation series. Using the BGI framework, approximate symmetries, first integrals, and approximate solutions have been constructed for a number of models involving ordinary and partial differential equations \cite{unal2000periodic, baikovjoining, bai2018approximate}. In Ref. \cite{tarayrah2021relationship}, the authors proved that to every point or local symmetry of an exact ODE of $n^{\rm th}$ order, there corresponds a BGI approximate local symmetry of the perturbed ODE of order at most $n-1$. A different approach to approximate symmetries, developed by Fushchich and Shtelen \cite{fushchich1989approximate}, combines a perturbation technique with the symmetry group method by expanding the dependent variables in a Taylor series in the small parameter. Using this method, approximate symmetries and approximate solutions have been found for some PDE models \cite{euler1992approximate,euler1994}. The BGI and Fushchich-Shtelen approaches have been compared and used to obtain approximate symmetries and approximate solutions for several PDE models \cite{wiltshire2006two,grebenev2007approximate}.

The paper is organized as follows. In Section \ref{sec:PAs}, we briefly overview the framework of Lie point and local symmetries of differential equations in comparison with approximate symmetry framework suggested by Baikov, Gazizov and Ibragimov for models involving a small parameter. In Section \ref{sec:Alg}, we find the BGI approximate point symmetries of the perturbed Boussinesq ODE, and we find the approximate local symmetries corresponding to the unstable point symmetries of the unperturbed equation, then we use some of the computed local symmetries to construct a general approximate solution of the Boussinesq ODE. In Section 4, we find an approximate integrating factor for a Benjamin-Bona-Mahony (BBM) perturbed ODE and use it to find an approximate solution of this equation. A brief discussion is offered in Section \ref{conclusion}.

\section{Lie groups of exact and approximate point and local symmetries}\label{sec:PAs}
A general system of $N$ differential equations is given by
\begin{equation}\label{gs1}
   F_0^{\sigma}(x,v,\partial v,\ldots, \partial^k v)=0,\quad k\geq 1, \quad \sigma=1,2,...,N,
\end{equation}
and its first-order perturbation in terms of a small parameter $\epsilon$ is written as:
\begin{equation}\label{eq:perturbed:general}
\barr
   F^{\sigma}(x,v,\partial v,\ldots, \partial^k v;\epsilon)=F_0^{\sigma}(x,v,\partial v,\ldots, \partial^k v)+\epsilon F_1^{\sigma}(x,v,\partial v,\ldots, \partial^k v)=O(\epsilon^2).
\earr
\end{equation}
Here $x = (x^1 ,x^2 ,...,x^n )$, $n\geq 1$, and $v=(v^1, v^2, ..., v^m)$, $m\geq 1$, denote respectively independent and dependent variables, and a symbol $\partial^q v$ denotes all $q^{\rm th}$-order derivatives of all components of $v$. A one-parameter Lie group of transformations
\begin{equation}\label{e25}
\begin{array}{l}
  (x^*)^i  = f^i (x,v;a) =x^i+a\xi_0^i(x,v)+ O(a^2) ,\quad i = 1,2, \ldots ,n, \\
  (v^*)^\mu  = g^\mu (x,v;a)= v^\mu +a \eta_0^\mu(x,v)+ O(a^2),\quad \mu = 1,2, \ldots ,m, \\
 \end{array}
\end{equation}
with the group parameter $a$, and the corresponding infinitesimal generator
\begin{equation}\label{e100}
X^0 = \xi_0^i (x,v)\dfrac{\partial }{{\partial x^i }}
+ \eta_0 ^\mu (x,v) \dfrac{\partial }{{\partial v^\mu }}
\end{equation}
is a point symmetry of the system \eqref{gs1} when for each $\sigma=1,2,...,N,$
\begin{equation}\label{det:eq:exact}
X^{0(k)} F_0^{\sigma}(x,v,\partial v,\ldots, \partial^k v)=0
\end{equation}
holds on solutions of \eqref{gs1} (e.g., \cite{olver2000applications, bluman2010applications}). The evolutionary (characteristic) form of the Lie group of transformations \eqref{e25} is the one-parameter family of transformations:
\begin{equation}\label{evol}
\begin{array}{l}
(x^*)^i  =x^i ,\quad i = 1,2, \ldots ,n, \\
(v^*)^\mu  = v^\mu +a [{\eta}_0^\mu(x,v)-v_i^\mu \xi_0^i(x,v)]+ O(a^2),\quad \mu = 1,2, \ldots ,m, \\
\end{array}
\end{equation}
with the evolutionary infinitesimal generator
\begin{equation}\label{evolg}
\hat{X}^0=[{\eta}_0^\mu(x,v)-v_i^\mu \xi_0^i(x,v)]\frac{\partial}{\partial v^\mu}.
\end{equation}
Higher-order local transformations generalize \eqref{evol} by allowing the infinitesimal components to depend on higher derivatives of $v$.

For equations \eqref{eq:perturbed:general} involving a small parameter $\epsilon$, we now define Baikov-Gazizov-Ibragimov (BGI) \cite{ibragimov1995crc,ibragimov2009approximate} approximate point symmetries . The BGI approximate symmetry generator for the perturbed equations \eqref{eq:perturbed:general} is given by
\beq\label{X:X0:X1}
X = X^0+\epsilon X^1= \left(\xi_0^i(x,v)+\epsilon \xi_1^i(x,v)\right)\dfrac{\partial}{\partial x^i}+\left(\eta_0^\mu(x,v)+\epsilon \eta_1^\mu(x,v)\right)\dfrac{\partial}{\partial v^\mu}.
\eeq
Similarly to exact higher-order local transformations, one can define local approximate BGI transformations, with generators in evolutionary form given by
\beq\label{X:X0:X1:higher}
\hat{X} = \hat{X}^0+\epsilon \hat{X}^1 = \left(\zeta_0^\mu[v]+\epsilon \zeta_1^\mu[v]\right)\dfrac{\partial}{\partial v^\mu}.
\eeq
The following theorem holds \cite{ibragimov1995crc}.
\begin{theorem}\label{1002}
Let the equations \eqref{eq:perturbed:general} be approximately invariant under an approximate group of BGI point transformations with the generator \eqref{X:X0:X1}
such that $\xi^0, \eta^0(x,v) \neq 0$. Then the infinitesimal operator \eqref{e100} is a generator of an exact symmetry group for the unperturbed equations \eqref{gs1}.
\end{theorem}
The determining equations to find the first-order BGI approximate symmetries of the system \eqref{eq:perturbed:general} are given by
\begin{equation}\label{detapprox}
\barr
    (X^{0^{(k)}}+\epsilon X^{1^{(k)}})(F_0^{\sigma}(x,v,\partial v,\ldots, \partial^k v)+\epsilon F_1^{\sigma}(x,v,\partial v,\ldots, \partial^k v))\bigg |_{F_0+\epsilon
    F_1=0}=O(\epsilon^2),\\[2ex]
    \sigma=1,...,N.
\earr
\end{equation}

 An exact point (or local) symmetry $X^0$ of the unperturbed equations \eqref{gs1} is called \textit{\textrm{stable}} if there exists a point (or local) generator $X^{1}$ such that
\eqref{X:X0:X1} (or respectively \eqref{X:X0:X1:higher}) is a BGI approximate symmetry of the perturbed equation \eqref{eq:perturbed:general}. If all symmetries of the equations \eqref{gs1} are \emph{stable}, the perturbed equations \eqref{eq:perturbed:general} are said to \emph{inherit} the symmetries of the unperturbed equations \cite{ibragimov1995crc}. In \cite{tarayrah2021relationship}, the authors proved that all exact symmetries of higher-order ODEs are stable and they correspond to higher-order approximate symmetries of the perturbed model.
\section{Approximate local symmetries and approximate solutions of the Boussinesq ODE} \label{sec:Alg}
In this section, we follow a systematic way developed in \cite{tarayrah2021relationship} to find higher-order BGI approximate symmetries of a second-order Boussinesq ODE that correspond to every point symmetry of the unperturbed equation and we use the approximate local symmetries to construct an approximate solution of the perturbed ODE.

 Consider a linear ODE
  \begin{equation}\label{unpert Boss}
    y''+y=0
  \end{equation}
and its perturbed version, the Boussinesq ODE
  \begin{equation}\label{pert Boss}
  y''+y-\epsilon \left(x+1+y^2\right)=0.
  \end{equation}
  The latter ODE can be obtained from the integration of a traveling wave reduction of the Boussinesq partial differential equation
\begin{equation}\label{BoussPDE}
  u_{tt}-u_{xx}+\epsilon (u^2)_{xx}-u_{xxxx}=0,\quad u=u(x,t),
\end{equation}
that was introduced in 1871 to describe the propagation of long waves in shallow water \cite{clarksonz1989}.

First, we seek exact point symmetries for \eqref{unpert Boss}. The exact point symmetry generator of the ODE \eqref{unpert Boss} has the form
\begin{equation}\label{exact Boss}
  X^0=\xi^0(x,y)\dfrac{\partial}{\partial x}+ \eta^0(x,y) \dfrac{\partial}{\partial y}.
\end{equation}
 Applying the determining equations \eqref{det:eq:exact}, one finds
 \begin{equation}\label{inf0}
   \barr
     \xi^0 = -C_1y\cos x+C_2y\sin x +C_3\sin 2x -C_4 \cos 2x +C_5,\\ [2 ex]
     \eta^0 = C_1 y^2\sin x+C_2 y^2\cos x+C_3y \cos 2x+C_4 y \sin 2x+C_6 \sin x +C_7 \cos x +C_8y.
  \earr
 \end{equation}
 Consequently, the ODE \eqref{unpert Boss} admits eight-dimensional Lie algebra of point symmetry generators, spanned by
 \begin{equation}\label{point symm Boss}
   \barr
    X_1^0 = y^2\sin x \dfrac{\partial}{\partial y}- y\cos x \dfrac{\partial}{\partial x},\quad X_2^0=  y^2\cos x \dfrac{\partial}{\partial y}+ y\sin x \dfrac{\partial}{\partial x},\quad X_3^0=  y\cos 2x \dfrac{\partial}{\partial y}+\sin 2x \dfrac{\partial}{\partial x},  \\
    X_4^0= y\sin 2x \dfrac{\partial}{\partial y}-\cos 2x \dfrac{\partial}{\partial x},\quad X_5^0= \dfrac{\partial}{\partial x},\quad X_6^0=\sin x\dfrac{\partial}{\partial y},\quad X_7^0=\cos x\dfrac{\partial}{\partial y},\quad X_8^0=y\dfrac{\partial}{\partial y}.
   \earr
 \end{equation}
Next, we find approximate point symmetries of the Boussinesq ODE \eqref{pert Boss}. Let
    \begin{equation}\label{}
      X = X^0+\epsilon X^1 = \left(\xi^0(x,y)+\epsilon \xi^1(x,y)\right)\dfrac{\partial}{\partial x}+ \left(\eta^0(x,y)+\epsilon \eta^1(x,y)\right) \dfrac{\partial}{\partial y}
    \end{equation}
denote the approximate BGI symmetry generator admitted by \eqref{pert Boss}, where $X^0$ is an exact symmetry generator \eqref{exact Boss} of the unperturbed ODE \eqref{unpert Boss}. The determining equation for approximate
symmetries \eqref{detapprox} yields
\begin{equation}\label{}
  \barr
  \xi^1 = -a_1y\cos x+a_2y\sin x +a_3\sin 2x -a_4 \cos 2x +a_5-\dfrac{4 C_6}{3} \cos x+\dfrac{4 C_7}{3} \sin x,\\
   \eta^1 = a_1 y^2\sin x+a_2 y^2\cos x+a_3y \cos 2x+a_4 y \sin 2x+a_6 \sin x +a_7 \cos x +a_8y\\+\dfrac{2 C_6}{3}y \sin x+\dfrac{2 C_7}{3}y \cos x+C_5.\\
  \earr
\end{equation}
The determining equation \eqref{detapprox} also provides some restrictions on the unperturbed symmetry components \eqref{inf0}: $C_1=C_2=C_3=C_4=C_8=0$. It follows that the point symmetries $X^0_j$, $j=1,\ldots, 4$ and $X_0^8$ in \eqref{point symm Boss} are unstable. The perturbed ODE \eqref{pert Boss} admits eight trivial approximate symmetries $X_j=\epsilon X^0_j,\, j=1,2,...,8,$ where $X^0_j$ are the exact point symmetries \eqref{point symm Boss} of the unperturbed ODE \eqref{unpert Boss}, and three nontrivial approximate point symmetries
   \begin{equation}\label{x7}
   \barr
   X_{9}=\dfrac{\partial}{\partial x}+\epsilon \dfrac{\partial}{\partial y},\quad X_{10}= \sin x\dfrac{\partial}{\partial y}+\epsilon \left(\dfrac{2}{3}y \sin x\dfrac{\partial}{\partial y}-\dfrac{4}{3} \cos x \dfrac{\partial}{\partial x}\right),\\ X_{11}= \cos x\dfrac{\partial}{\partial y}+\epsilon \left(\dfrac{2}{3}y \cos x\dfrac{\partial}{\partial y}+\dfrac{4}{3} \sin x \dfrac{\partial}{\partial x}\right)
   \earr
   \end{equation}
   corresponding respectively to the stable point symmetries $X_0^5$, $X_0^6$ and $X_0^7$ of the unperturbed equation \eqref{unpert Boss}.

   As guaranteed in \cite{tarayrah2021relationship}, the unstable point symmetries of the unperturbed equation \eqref{unpert Boss} correspond to first-order approximate local symmetries of the Boussinesq ODE \eqref{pert Boss}. The corresponding approximate local symmetries for each unstable point symmetry of \eqref{unpert Boss} takes the form
    \[
  \hat{X} =\left(\zeta^0+\epsilon \zeta^1(x,y,y')\right)\dfrac{\partial}{\partial y},
\]
where $\zeta^0{\partial}/{\partial y}$ is the evolutionary form of the point symmetry generator \eqref{exact Boss} of \eqref{unpert Boss}. For the unstable point symmetry $X_0^1$, $\zeta^0= y^2\sin x+ yy'\cos x$. The corresponding $\zeta^1$ is any solution of the linear PDE
\[
 \left(D^2 \zeta^{1}+\zeta^{1}\right)\bigg|_{y''=-y}=\left(\left(3y^2+3x+3\right)y'+y\right)\cos x-2y^3\sin x,
\]
which has a particular solution given by
\[
 \zeta^{1}(x,y,y')=\left(\dfrac{2}{3}(y')^3+(y^2-x-1)y'-y\right)\cos x-2(x+1)y \sin x.
\]
One obtains
\[
  \hat{X}^1 =\left[y^2\sin x+ yy'\cos x+\epsilon \left(\left(\dfrac{2}{3}(y')^3+(y^2-x-1)y'-y\right)\cos x-2(x+1)y \sin x\right)\right]\dfrac{\partial}{\partial y}
\]
as a first-order approximate symmetry for the perturbed ODE \eqref{pert Boss} corresponds to the {\textrm{unstable}} point symmetry $X_0^1$. Similarly, one can obtain the higher-order approximate symmetries that correspond to the unstable point symmetries $X_0^2, X_0^3, X_0^4$ and $X_0^8$ in \eqref{point symm Boss}. They are given by
\begin{subequations}\label{symm:bouss:approx:3rd}
\beq
   \hat{X}^2 = \left[y^2\cos x- yy'\sin x+\epsilon \left(\left(-\dfrac{2}{3}(y')^3+(-y^2+x+1)y'+y\right)\sin x-2(x+1)y \cos x\right)\right]\dfrac{\partial}{\partial y},
\eeq
\beq
  \hat{ X}^3 = \left[ y\cos 2x-y'\sin 2x+\epsilon\left(\left(\dfrac{y^2}{3}-2(y')^2-x-1\right)\cos 2x+\left(-\dfrac{8}{3} yy'+1\right)\sin 2x\right)\right] \dfrac{\partial}{\partial y},
\eeq
\begin{equation}\label{}
 \hat{X}^4 =\left[ y\sin 2x+y'\cos 2x+\epsilon\left(\left(\dfrac{y^2}{3}-2(y')^2-x-1\right)\sin 2x+\left(\dfrac{8}{3} yy'-1\right)\cos 2x\right)\right] \dfrac{\partial}{\partial y},
\end{equation}
\begin{equation}\label{y-symm}
 \hat{X}^8 = \left[ y+\epsilon\left(\dfrac{2(y')^2}{3}+\dfrac{y^2}{3}-x-1\right)\right] \dfrac{\partial}{\partial y}.
\end{equation}
 \end{subequations}
 Now we proceed to find approximate solution of the ODE \eqref{pert Boss} using first-order approximate symmetries admitted by the Boussinesq ODE \eqref{pert Boss}. The fundamental solution of the unperturbed equation \eqref{unpert Boss} is
 \begin{equation}\label{IF33}
     y(x)=C_1\sin x+C_2 \cos x.
 \end{equation}
 The solution \eqref{IF33} is invariant under the group generated by
 \begin{equation}\label{IF34}
   X_8^0-C_1 X_6^0 - C_2 X_7^0 =\left(y-C_1\sin x-C_2 \cos x\right)\dfrac{\partial}{\partial y},
 \end{equation}
 where $ X_j^0,\,j=6,7,8$ are the point symmetries \eqref{point symm Boss} for the unperturbed ODE \eqref{unpert Boss}. $X_8^0$ is {\textrm{unstable}} point symmetry, however it corresponds to first-order approximate symmetry $\hat{X}^8$ in \eqref{y-symm}. $X_6^0$ and $X_7^0$ are {\textrm{stable}} as point symmetries, the corresponding approximate symmetries are given respectively by
 \[\hat{X}^6=\left(\sin x+\epsilon \left(\dfrac{2}{3}y \sin x+\dfrac{4}{3} y'\cos x \right) \right)\dfrac{\partial}{\partial y},\quad \hat{X}^7=\left(\cos x+\epsilon \left(\dfrac{2}{3}y \cos x-\dfrac{4}{3} y'\sin x \right) \right)\dfrac{\partial}{\partial y}.\]
 The approximately invariant solution under $\hat{X}^8-C_1\hat{X}^6-C_2\hat{X}^7$ is defined by
   \begin{equation}\label{IF 34}
     y-C_1x-C_2 \sin x-C_3 \cos x-C_4+\epsilon g(x,y,y')=O(\epsilon^2),
   \end{equation}
    where $g$ is given by
    \begin{equation}\label{}
      g=\dfrac{2(y')^2}{3}+\dfrac{y^2}{3}-x-1-C_1\left(\dfrac{2}{3}y \sin x+\dfrac{4}{3} y'\cos x\right)-C_2 \left(\dfrac{2}{3}y \cos x-\dfrac{4}{3} y'\sin x\right).
    \end{equation}
     Substituting $y(x;\epsilon)=y_0(x)+\epsilon y_1(x)+o(\epsilon)$ into the equation \eqref{IF 34}, one finds the following approximate solution of the Boussinesq ODE \eqref{pert Boss}
     \begin{equation}\label{IF35}
       y(x;\epsilon)=C_1 \sin x + C_2 \cos x +\epsilon \left(x+1+\dfrac{(C_1^2-C_2^2)\cos^2 x-C_1C_2 \sin 2x+C_1^2+2C_2^2}{3}\right).
     \end{equation}

\begin{figure}[H]
\centering
	\includegraphics[width=0.7\textwidth]{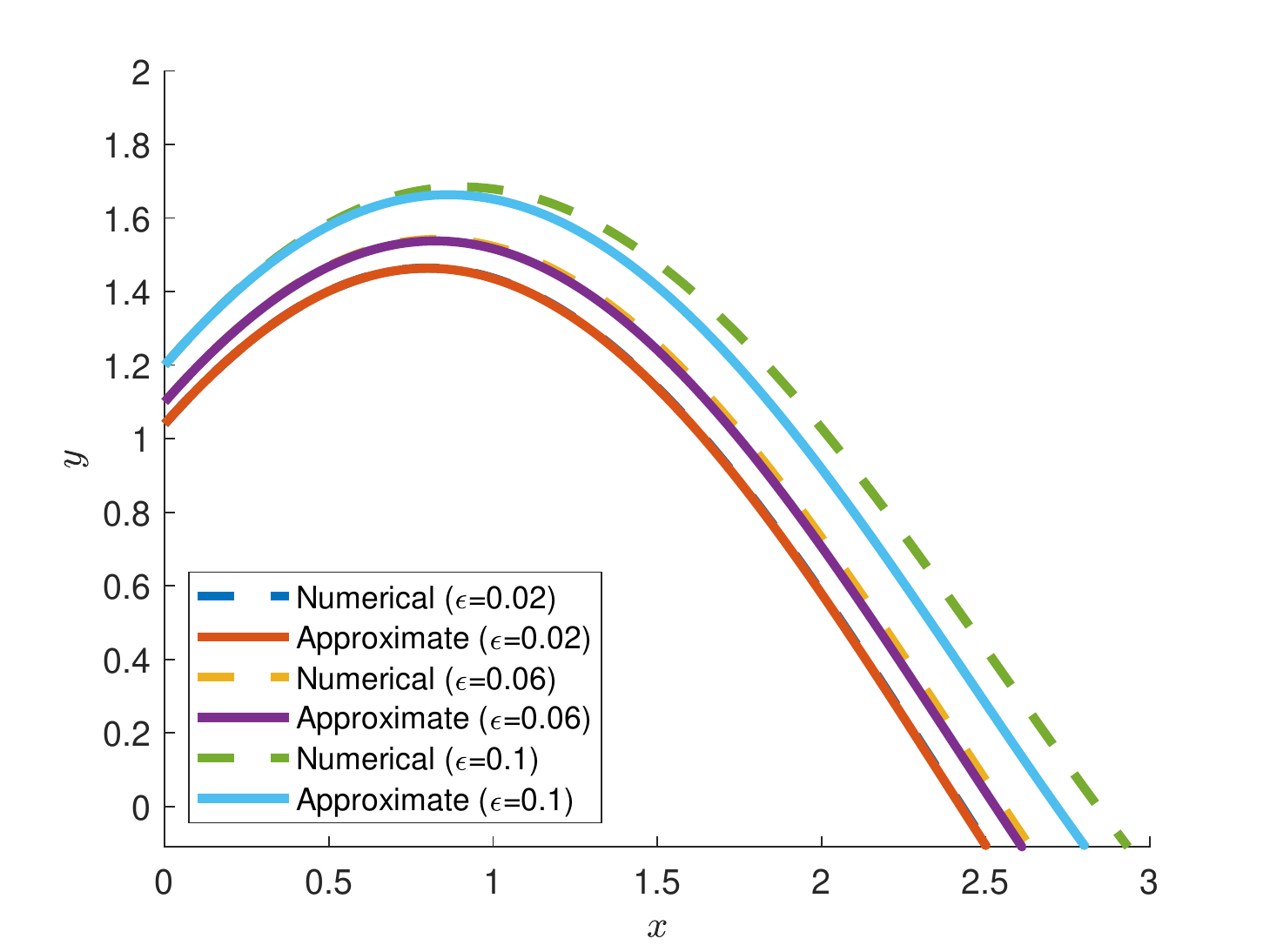}
	\caption{ The approximate solutions \eqref{Bouss particular approx} of the perturbed equation \eqref{pert Boss} with initial conditions \eqref{ICs2} (solid lines) vs.~numerical solutions of the same initial value problem (dashed lines) for \,$\epsilon=0.02, 0.06$ and $0.1$ (bottom to top).}\label{fig_sec7}
\end{figure}

The Boussinesq ODE \eqref{pert Boss} with initial conditions
 \begin{equation}\label{ICs2}
  y(0)=1+2 \epsilon,\quad y'(0)=1+\dfrac{\epsilon}{3}
\end{equation}
has a particular approximate solution given by
\begin{equation}\label{Bouss particular approx}
  y(x;\epsilon)=\sin x +\cos x +\epsilon \left(x+2-\dfrac{\sin 2x}{3}\right).
\end{equation}
Now, we compute a numerical solution for the perturbed ODE \eqref{pert Boss} with the same initial conditions \eqref{ICs2} using the \verb"Matlab" solver (ode45). It is an embedded method from the fourth- and fifth-order Runge-Kutta methods \cite{butcher2008numerical} with a variable time step for efficient computation based on an algorithm of Dormand and Prince \cite{dormand1980family}. In Figure \ref{fig_sec7}, graphs of the approximate solution \eqref{Bouss particular approx} vs.~the numerical solution for the initial conditions \eqref{ICs2} are shown for different values of $\epsilon$, showing good agreement.
\section{Approximate solution of Benjamin-Bona-Mahony (BBM) ODE using approximate integrating factors}
The Benjamin-Bona-Mahony (BBM) equation is the partial differential equation
\begin{equation}\label{bbm-pde}
  u_t+u_x-u_{xxt}+\dfrac{3}{2}\epsilon uu_x=0,\quad u=u(x,t).
\end{equation}
This equation was studied in \cite{benjamin1972model} as an alternative of the Korteweg–de Vries equation (KdV equation) for modeling long surface gravity waves of small amplitude $(\epsilon)$ propagating unidirectionally in $1+1$ dimensions. Using a travelling wave solution $u(x,t)=y(z)=y(x-ct)$,
equation \eqref{bbm-pde} reduces to an ODE:
\begin{equation}\label{bbm-ode}
  y'''+\dfrac{1-c}{c}y'+\dfrac{3}{2c}\epsilon yy'=0,
\end{equation}
The BBM ODE \eqref{bbm-ode} can be simplified to
\begin{equation}\label{bbm-ode1}
  y''+\dfrac{1-c}{c}y+\dfrac{3}{4c}\epsilon y^2=k,
\end{equation}
where $k$ is a constant of integration. By taking $k={3\epsilon}/{4c}$, the ODE \eqref{bbm-ode1} becomes
\begin{equation}\label{bbm-ode2}
   y''=\dfrac{c-1}{c}y-\dfrac{3}{4c}\epsilon (y^2-1).
\end{equation}
The approximate integrating factor for the second-order perturbed ODE
\begin{equation}\label{IF18}
  y''=f_0(x,y,y')+\epsilon f_1(x,y,y')
\end{equation}
has the form $\mu(x,y,y';\epsilon)=\mu_0(x,y,y')+\epsilon \mu_1(x,y,y')$, where the components $\mu_0$, $\mu_1$ satisfy the following equations \cite{tarayrah2021relationship}:
\begin{subequations}\label{eq:th51}
\beq
y'\mu_{0_{yy'}}+\mu_{0_{xy'}}+2\mu_{0_{y}}+(\mu_0 f_0)_{y'y'} = 0,\label{IF20}
\eeq
\beq
y'^2\mu_{0_{yy}}+2y'\mu_{0_{xy}}+\mu_{0_{xx}}+y'(\mu_0 f_0)_{yy'}+(\mu_0 f_0)_{xy'}- (\mu_0 f_0)_{y}= 0,\label{IF21}
\eeq
\beq
 y'\mu_{1_{yy'}}+\mu_{1_{xy'}}+2\mu_{1_{y}}+(\mu_1 f_0)_{y'y'}+(\mu_0 f_1)_{y'y'} = 0,\label{IF22}
\eeq
\beq \label{IF23}
\barr
y'^2\mu_{1_{yy}}+2y'\mu_{1_{xy}}+\mu_{1_{xx}}+y'(\mu_1 f_0)_{yy'}+(\mu_1 f_0)_{xy'}-(\mu_1 f_0)_{y}-(\mu_0 f_1)_{y} \\[2ex]
\qquad +y'(\mu_0 f_1)_{yy'}+(\mu_0 f_1)_{xy'} = 0.
\earr
\eeq
\end{subequations}
A possible approximate integrating factor for the perturbed ODE \eqref{bbm-ode2} is $\mu=(1+\epsilon)y'$. Multiplying this integrating factor by \eqref{bbm-ode2} and then integrating the resulting equation, one finds an approximate first integral:
  \begin{equation}\label{}
     D\left( y'^2+\dfrac{1-c}{c}y^2+\epsilon \left(y'^2+\dfrac{1-c}{c}y^2-\dfrac{3y}{2c}+\dfrac{y^3}{2c}\right)\right)=O(\epsilon^2).
  \end{equation}
  Consequently, the perturbed BBM ODE \eqref{bbm-ode2} is reduced to the first-order ODE
  \begin{equation}\label{IF30}
   y'^2+\dfrac{1-c}{c}y^2+\epsilon \left(y'^2+\dfrac{1-c}{c}y^2-\dfrac{3y}{2c}+\dfrac{y^3}{2c}\right)=C_1+O(\epsilon^2).
  \end{equation}
  Substituting $y(z;\epsilon)=y_0(z)+\epsilon y_1(z)+O(\epsilon^2)$ into the ODE \eqref{IF30} leads to the system of ODEs
\begin{equation}\label{sys-sec4}
  \barr
    (y_0')^2+\dfrac{1-c}{c}y_0^2 = C_1, \\[1ex]
    2y_0'y_1'+\dfrac{2(1-c)}{c}y_0y_1+(y_0')^2+\dfrac{1-c}{c}y_0^2-\dfrac{3y_0}{2c}+\dfrac{y_0^3}{2c} = 0.
\earr
\end{equation}
 The solution of the above system is bounded if $0<c<1$. Hence when $0<c<1$, the solution of the system \eqref{sys-sec4} is given by
 \begin{equation}\label{sol-sec4}
   \barr
   y_0(z)=\left(C_1c^2+1\right)\sin\left(\sqrt{\dfrac{1-c}{c}}(z-C_2)\right),\\[1ex]
   y_1(z)=\dfrac{1}{16c(1-c)^2}\biggl[ 4\sqrt{c(1-c)}\left(C_1c^{5/2}-C_1c^{3/2}+c^{1/2}-1\right)\sin\left(\sqrt{\dfrac{1-c}{c}}(z-C_2)\right)\\-(C_1c^2+1)\cos^2\left(\sqrt{\dfrac{1-c}{c}}(z-C_2)\right)+16C_3c(1-c)^2\cos\left(\sqrt{\dfrac{1-c}{c}}(z-C_2)\right)\\-C_1^2c^4-2C_1c^2-12c^2+12c-1\biggl].
 \earr
 \end{equation}
Consequently, a general approximate solution for the BBM ODE \eqref{bbm-ode} is given by
 \begin{equation}\label{BBM_approx_soln}
 \barr
   y(z;\epsilon)=\left(C_1c^2+1\right)\sin\left(\sqrt{\dfrac{1-c}{c}}(z-C_2)\right)\\[2ex]
   +\dfrac{\epsilon}{16c(1-c)^2}\biggl[ 4\sqrt{c(1-c)}\left(C_1c^{5/2}-C_1c^{3/2}+c^{1/2}-1\right)\sin\left(\sqrt{\dfrac{1-c}{c}}(z-C_2)\right)\\-(C_1c^2+1)\cos^2\left(\sqrt{\dfrac{1-c}{c}}(z-C_2)\right)+16C_3c(1-c)^2\cos\left(\sqrt{\dfrac{1-c}{c}}(z-C_2)\right)\\-C_1^2c^4-2C_1c^2-12c^2+12c-1\biggl].
 \earr
 \end{equation}
Hence the general approximate solution of the BBM PDE \eqref{bbm-pde} is $u(x,t)=y(x-ct;\epsilon)$.

We could not find harmonic-type solutions of BBM PDE in literature. The exact solutions of different forms of BBM equation is given in terms of Jacobi elliptic functions $\textrm{cn}(v,k), \textrm{sn}(v,k), 0\leq k \leq 1.$ When $k\rightarrow 1,$ one obtains the solitary wave solutions of BBM equation (see, e.g.\cite{fu2003jefe,wazwaz2006new,chen2007cnoidal}). For the BBM model \eqref{bbm-pde}, the explicit cnoidal wave solutions are given by
\begin{eqnarray}\label{cn}
  u(x,t) &=& a^2{\textrm{cn}}^2(B_1(x-ct),k)+H_1,  \\
  u(x,t) &=& a^2{\textrm{sn}}^2(B_2(x-ct),k)+H_2,\label{sn}
\end{eqnarray}
where $B_1$, $B_2$, $H_1$ and $H_2$ are given in terms of $a, \epsilon, c$ and $k$:
\[B_1=\dfrac{a\sqrt{2\epsilon c }}{4ck},\quad H_1=\dfrac{2ck^2+\epsilon a^2-2 \epsilon a^2k^2-2k^2}{3 \epsilon k^2},\]
\[B_2=\dfrac{a\sqrt{-2\epsilon c }}{4ck},\quad H_2=\dfrac{2ck^2-\epsilon a^2- \epsilon a^2k^2-2k^2}{3 \epsilon k^2}.\]
When $k\rightarrow 1$, \eqref{cn} reduces to a solitary wave solution:
\begin{equation}\label{}
  u(x,t)=a^2 {\sech}^{2} \left(\frac{a\sqrt{2\epsilon c  }\,(x-ct)}{4c}\right)+\dfrac{2c-\epsilon a^2-2}{3 \epsilon}.
\end{equation}
When $c<0$ and $k\rightarrow 1$, the solution \eqref{sn} simplifies to a left-travelling kink wave solution:
\begin{equation}\label{}
 u(x,t)=a^2 {\tanh}^{2} \left(\frac{a\sqrt{-2\epsilon c }\,(x-ct)}{4c}\right)+\dfrac{2c-2\epsilon a^2-2}{3 \epsilon}.
\end{equation}
Note that harmonic-type solutions like \eqref{BBM_approx_soln} do not follow from \eqref{cn} and \eqref{sn} as $k\rightarrow 0^{+}$. With the initial conditions
\begin{equation}\label{ICS3}
  y(0)=-\dfrac{\epsilon}{16},\quad y'(0)=\dfrac{5}{4}-\dfrac{5\epsilon}{8}, \quad y''(0)= \dfrac{25\epsilon}{16},
\end{equation}
the BBM ODE \eqref{bbm-ode} has a particular solution given by
\begin{equation}\label{approx_part_BBM}
  y(x;\epsilon)=\dfrac{5}{4}\sin x+\epsilon \left(\dfrac{23-25\cos^2 x-20 \sin x}{32}\right).
\end{equation}
Now, we compute a numerical solution for the BBM ODE \eqref{bbm-ode} with the initial conditions \eqref{ICS3} using the \verb"Matlab" solver (ode45). In Figure \ref{fig2} below, graphs with the solid lines refer to the numerical solution of \eqref{bbm-ode} while the long
dash styles represent the graphs of the approximate solution \eqref{approx_part_BBM} of the perturbed equation \eqref{bbm-ode} for different values of the small parameter $\epsilon$.
\begin{figure}[H]
\centering
	\includegraphics[width=0.7\textwidth]{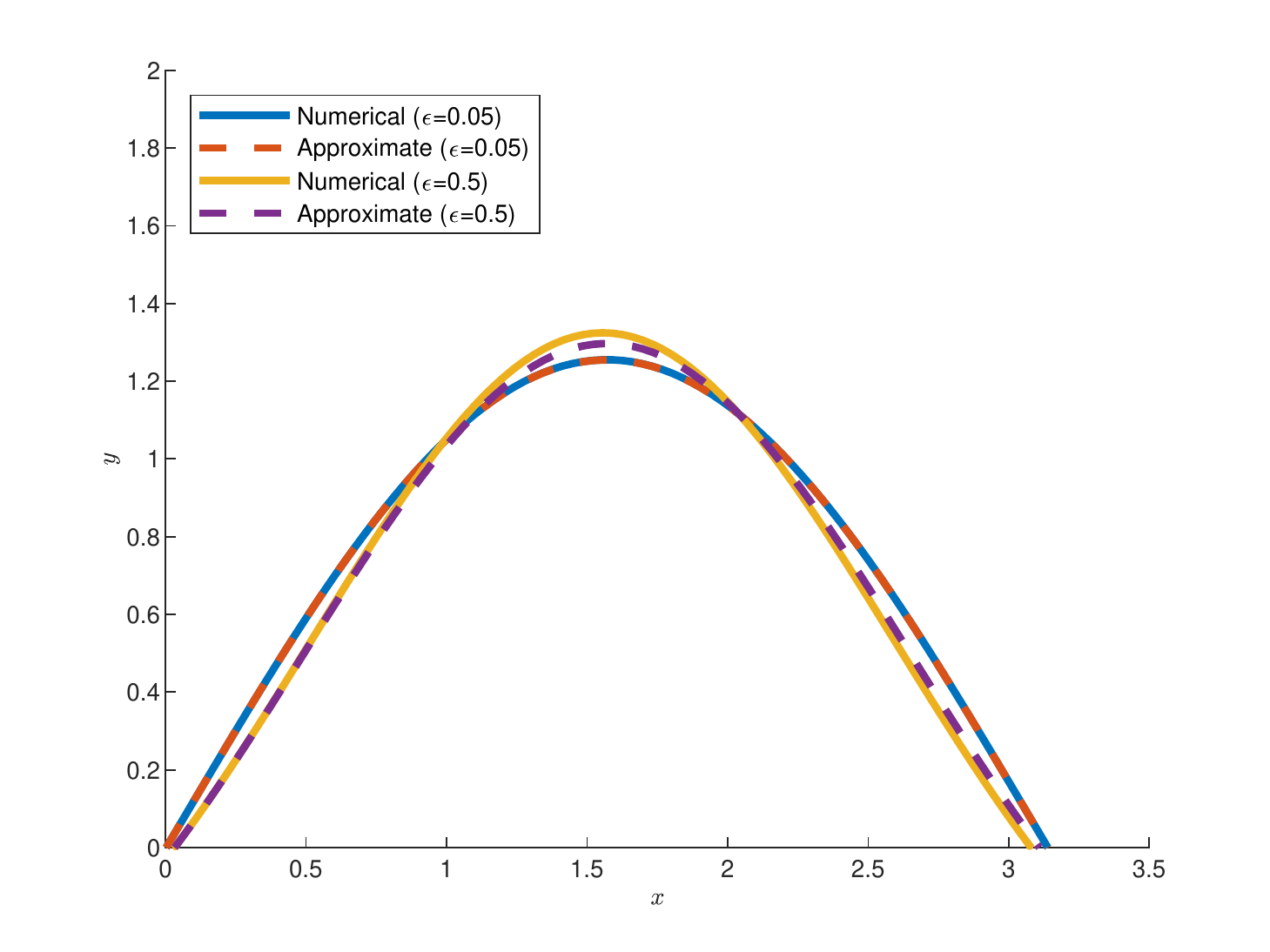}
	\caption{The approximate solutions \eqref{approx_part_BBM} of the perturbed equation \eqref{bbm-ode} vs. the numerical
solutions of \eqref{bbm-ode} for $\epsilon= 0.05, 0.5$ with initial conditions \eqref{ICS3}.}\label{fig2}
\end{figure}

\section{Discussion}  \label{conclusion}

In this paper, we considered approximate local symmetries (in the BGI framework \cite{baikov1989, baikov1991, baikov1993}) and approximate integrating factors of ordinary differential equations involving a small parameter ($\epsilon$). We found approximate local symmetries of the Boussinesq perturbed ODE \eqref{pert Boss} that correspond to every unstable point symmetry of the unperturbed equation \eqref{unpert Boss}. These approximate local symmetries were used to construct a general approximate solution of the ODE \eqref{pert Boss}. A particular approximate solution for the Boussinesq ODE \eqref{pert Boss} with initial conditions \eqref{ICs2} was found and compared to a numerical solution of \eqref{pert Boss} with the same initial conditions. A sine-type approximate solution \eqref{BBM_approx_soln} of the BBM perturbed ODE \eqref{bbm-ode} was constructed using an approximate integrating factor for \eqref{bbm-ode} and this solution could not be obtained from the cnoidal wave solutions \eqref{cn} and \eqref{sn} of the PDE \eqref{bbm-pde}.

\section*{Acknowledgments }
The author is grateful to Prof. Alexei Cheviakov for the encouragement and discussion of the results.


{\footnotesize
\bibliography{References_3}
\bibliographystyle{ieeetr}}
\end{document}